\documentclass[11pt]{amsart}
\usepackage[hmargin=2.5cm,vmargin=2.5cm]{geometry}
\usepackage{amsfonts}
\usepackage{amsthm}
\usepackage[english]{algorithm2e}
\usepackage{amsmath}
\usepackage{amssymb}
\usepackage[T1]{fontenc}
\newtheorem{Theorem}{Theorem}[section]
\newtheorem{Definition}[Theorem]{Definition}
\newtheorem{Remark}[Theorem]{Remark}

\newtheorem{Lemma}[Theorem]{Lemma}

\title{Another super-identity equivalent to the Hom-Malcev super-identity}
\author[Sylvain Attan, Donatien Gaparayi and A. Nourou Issa]
       {Sylvain Attan $^1$,  Donatien Gaparayi $^2$ and A. Nourou Issa $^3$
       \\\\
          $^{1,3}$ D\'epartement de Math\'ematiques\\
       Universit\'e d'Abomey-calavi\\
       01 BP 4521, Cotonou 01, B\'enin\\
       $^2$ Ecole Normale Sup\'erieure (E.N.S),\\
        BP 6983 Bujumbura, Burundi\\
       $^1$ syltane2010@yahoo.fr, $^2$ gapadona@yahoo.fr, $^3$ woraniss@yahoo.fr
       }
\begin{document}
\maketitle

\begin{abstract}
In a Hom-superalgebra a super-identity, equivalent to the Hom-Malcev \\super-identity, is found.
\end{abstract}

\section{Introduction}
The theory of Hom-algebras originated from Hom-Lie algebras introduced by J.T. Hartwig, D. Larsson, and S.D. Silvestrov in\cite{Hartwig} in the study of
quasi-deformations of Lie algebras of vector fields, including q-deformations of Witt algebras and Virasoro algebras. As genera-\\lization of Hom-Lie 
algebras, Hom-Malcev algebras are introduced \cite{YAU2} to study Hom-alternative algebras \cite{MAK2}. A $\mathbb{Z}_2$-graded version of Hom-Lie 
(rep. Hom-Malcev) algebra called Hom-Lie (resp. Hom-Malcev) superalgebra are introduced \cite{Ammar} (resp. \cite{Nan}, \cite{Abdaoui}) as generalization of Lie (rep. Malcev) superalgebras. 
Recall that a Malcev superalgebra is a non-associative superalgebra $A$ with a super skewsymmetric multiplication
$\cdot$ (i.e, $xy =-(−1)^{\bar{x}\bar{y}}yx$) such that the Malcev super-identity
\begin{eqnarray}
 && 2tJ_A (x, y, z)=J_A (t, x, yz) + (-1)^{\bar{x}(\bar{y}+\bar{z})}J_A (t, y,zx)\nonumber\\
&&(-1)^{\bar{z}(\bar{x}+\bar{y})}J_A (t, z,xy)\label{a}
\end{eqnarray}
is satisfied for all homogeneous elements $x, y, z, t$ in the superspace $A,$ where
\begin{eqnarray}
 J_A(x, y, z)=((xy)z+(-1)^{\bar{x}(\bar{y}+\bar{z})}(yz)x+(−1)^{\bar{z}(\bar{x}+\bar{y})}(zx)y
\end{eqnarray}
is the super-Jacobian. In particular, Lie superalgebras are examples of Malcev superalgebras. Malcev
superalgebras play an important role in the geometry of smooth loops.\\

Some twisting of the Malcev super-identity (\ref{a}) along any even algebra self-map $\alpha$ of $A$ gives rise to the notion of a
 Hom-Malcev superalgebra $(A,[,], \alpha)$ (\cite{Abdaoui}; see definitions in section 2). Properties and constructions of Hom-Malcev superalgebras, as
 well as the relationships between these Hom-superalgebras and Hom-alternative or Hom-Jordan superalgebras are investigated in \cite{Abdaoui}. 
In particular, it is shown that a Malcev superalgebra can be twisted into a Hom-Malcev superalgebra and that Hom-alternative superalgebras 
are Hom-Malcev super-admissible.
In \cite{Abdaoui}, as for Malcev algebras (see \cite{Sagle}, \cite{LY1}), Hom-Malcev algebras \cite{YAU2},\cite{Issa1} and Malcev superalgebra\cite{Elduque2}, equivalent defining identities 
of a Hom-Malcev superalgebra are given. In this note, we mention another super identity in a Hom-Malcev superalgebra that is equivalent to
 the ones found in \cite{Abdaoui}. Specifically, we shall prove the following\\
\begin{Theorem} \rm{
  Let $(M=M_0\oplus M_1, \cdot, \alpha)$ be a Hom-Malcev superalgebra. Then the super identity
\begin{eqnarray}
 \tilde{J}_M(wx,\alpha(y),\alpha(z))&=&\alpha^2(w)\tilde{J}_M(x,y,z)+
(-1)^{\bar{x}(\bar{y}+\bar{z})}\tilde{J}_M(w,y,z)\alpha^2(x)\nonumber\\
   &&-2(-1)^{(\bar{y}+\bar{z})(\bar{x}+\bar{w})}\tilde{J}_M(yz,\alpha(w),\alpha(x))\label{s1}
\end{eqnarray}
holds for all homogeneous elements $w, x, y, z$ in $M.$}
\end{Theorem}
Moreover, in any super anticommutative Hom-superalgebra $(M=M_0\oplus M_1, \cdot, \alpha)$, the super identity (\ref{s1})
is equivalent to the Hom-Malcev super identity
\begin{eqnarray}
 && 2\alpha^2(t)\tilde{J}_M (x, y, z)=\tilde{J}_M (\alpha(t), \alpha(x), yz) + (-1)^{\bar{x}(\bar{y}+\bar{z})}\tilde{J}_M (\alpha(t), \alpha(y), zx)\nonumber\\
&&(-1)^{\bar{z}(\bar{x}+\bar{y})}\tilde{J}_M (\alpha(t), \alpha(z), xy)\label{b}
\end{eqnarray}
or equivalently
\begin{eqnarray}
 &&\tilde{J}_M(\alpha(y),\alpha(z),wx)+(-1)^{\bar{y}\bar{z}+\bar{w}(\bar{y}+\bar{z})}\tilde{J}_M(\alpha(w),\alpha(z),yx)\nonumber\\
&=&+(-1)^{\bar{w}\bar{x}}\tilde{J}_M(y,z,x)\alpha^2(w)+(-1)^{\bar{y}(\bar{z}+\bar{w}+\bar{x})+\bar{z}\bar{w}}\tilde{J}_M(w,z,x)\alpha^2(y)\label{c}
\end{eqnarray}
for all homogeneous elements $x, y, z, t$ in the superspace $M.$ (See \cite{Abdaoui} for equivalence between (\ref{b}), (\ref{c}) and another super identiry).\\
\\
Observe that when $\alpha= Id$ (the identity map) in (\ref{b}), then (\ref{b}) is (\ref{a})
i.e. the Hom-Malcev superalgebra $(M=M_0\oplus M_1, \cdot, \alpha)$ reduces to the Malcev superalgebra $(M=M_0\oplus M_1, \cdot)$
(see \cite{Abdaoui}). Also if all the homogeneous elements are in $M_0,$ (i.e are even) then (\ref{s1}) becomes 
\begin{eqnarray}
 \tilde{J}_M(wx,\alpha(y),\alpha(z))&=&\alpha^2(w)\tilde{J}_M(x,y,z)+
\tilde{J}_M(w,y,z)\alpha^2(x)\nonumber\\
   &&-2\tilde{J}_M(yz,\alpha(w),\alpha(x))\nonumber
\end{eqnarray}
which is an equivalent identity to Hom-Malcev identity\cite{YAU2} found in \cite{Issa1} i.e $(M=M_0\oplus M_1, \cdot, \alpha)$ reduces to Hom-Malcev 
algebras $(M_0,\cdot, \alpha).$\\
 In section 2 some instrumental lemmas are proved. Some results in these
lemmas are a kind of the $\mathbb{Z}_2$-graded version of similar results found in \cite{Issa1}
in case of Hom-Malcev algebras. The section 3 is devoted to the proof of the theorem.

\section{On Hom-superalgebras}
Throughout this paper, $\mathbb{K}$ is an algebraically closed field of characteristic 0 and $M$ is a linear super-espace
over $\mathbb{K}.$
In this section we recall useful notions on Hom-Lie superalgebras \cite{Ammar}
as well as the one of a Hom-Malcev superalgebra \cite{Abdaoui}. The main result
of this section (Lemma \ref{sis4}) proves that the $\mathbb{Z}_2$-graded version of the identity (2.4)
of Lemma 2.7 \cite{Issa1} which is the Hom-version of the identity (6) of \cite{Kleinf} holds in any Hom-Malcev superalgebra.\\
Now let $M$ be a linear superespace over $\mathbb{K}$ that is a $\mathbb{Z}_2$-graded linear space with a direct sum 
$M=M_0\oplus M_1.$ 
The element of $M_j,$ $j\in \mathbb{Z}_2,$ are said to be homogeneous of parity $j.$ The parity of a homogeneous element 
$x$ is denoted by $\bar{x}.$ In the sequel, we will denote by $\mathcal{H}(M)$ the set of all homogeneous elements of $M.$
\begin{Definition}\cite{Abdaoui}
 A Hom-superalgebra is a triple $(M,\mu, \alpha)$ in which $M$ is a $\mathbb{K}$-super-module, $\mu :
M^{\otimes 2}\rightarrow M$ is an even bilinear map, and $\alpha: M\rightarrow M$ is an even linear map such that 
$\alpha\circ\mu=\mu\circ\alpha^{\otimes 2}$ (multiplicativity)
\end{Definition}
\begin{Definition}\cite{Ammar} \rm {
  Let $(M,\cdot, \alpha)$ be a Hom-superalgebra.\\
(1) The Hom-superJacobien is the trilinear map
 $\tilde{J}_M$ defined as
\begin{eqnarray}
\tilde{J}_M(x,y,z)&=&(xy)\alpha(z)+(-1)^{\bar{x}(\bar{y}+\bar{z})}(yz)\alpha(x)+(-1)^{\bar{z}(\bar{x}+\bar{y})}(zx)\alpha(y)
\end{eqnarray}
 for all $x,y,z \in \mathcal{H}(M).$\\
(2) $(M,\cdot, \alpha)$ is called Hom-Lie superalgebra if $"\cdot"$ is super anticommutative (i.e $xy=-(-1)^{\bar{x}\bar{y}}yx$)
and $\tilde{J}_M(x,y,z)=0$ 
for all $x,y,z \in \mathcal{H}(M).$}
\end{Definition} 
\begin{Definition}\cite{Abdaoui}
 A Hom-Malcev superalgebra is a super anticommutative Hom-superalgebra 
$(M,\cdot, \alpha)$ such that 
\begin{eqnarray}
 && 2\alpha^2(t)\tilde{J}_M (x, y, z)=\tilde{J}_M (\alpha(t), \alpha(x), yz) + (−1)^{\bar{x}(\bar{y}+\bar{z})}\tilde{J}_M (\alpha(t), \alpha(y), zx)\nonumber\\
&&(−1)^{\bar{z}(\bar{x}+\bar{y})}\tilde{J}_M (\alpha(t), \alpha(z), xy)\nonumber\\
\end{eqnarray}
(see (\ref{b}) above) for all $x,y,z \in \mathcal{H}(M).$
\end{Definition}
\begin{Remark} \rm {
 Observe that any Hom-Lie superalgebra is a Hom-Malcev superalgebra since the Hom-Lie superidentity implies 
the one of Hom-Malcev}
\end{Remark}
We have  the following 
\begin{Lemma}\label{sis1} \rm{
 In any anticommutative Hom-superalgebra $(A=A_0\oplus A_1, \cdot, \alpha),$ the \\following holds\\
(i) $\tilde{J}_A(x,y,z)$ is super skew-symmetric in its three variables that is\\
$\tilde{J}_A(x,y,z)=-(-1)^{\bar{x}\bar{y}}\tilde{J}_A(y,x,z)=-(-1)^{\bar{y}\bar{z}}\tilde{J}_A(x,z,y)=
-(-1)^{\bar{x}(\bar{y}+\bar{z})+\bar{y}\bar{z}}\tilde{J}_A(z,y,x)$
 \begin{eqnarray}
 \mbox{(ii)}&&\alpha^2(w)\tilde{J}_A(x,y,z)-(-1)^{\bar{w}(\bar{x}+\bar{y}+\bar{z})}\alpha^2(x)\tilde{J}_A(y,z,w) +
(-1)^{(\bar{y}+\bar{z})(\bar{w}+\bar{x})} \alpha^2(y)\tilde{J}_A(z,w, x)\nonumber\\
&&-(-1)^{\bar{z}(\bar{x}+\bar{y}+\bar{w})}\alpha^2(z)\tilde{J}_A(w,x,y)\nonumber\\
&=& \tilde{J}_A(wx,\alpha(y),\alpha(z))+(-1)^{(\bar{y}+(\bar{z})(\bar{x}+\bar{w})}\tilde{J}_A(yz,\alpha(w),\alpha(x))
+(-1)^{\bar{x}(\bar{y}+\bar{z})}\tilde{J}_A(wy,\alpha(z),\alpha(x))\nonumber\\
&&+(-1)^{\bar{z}(\bar{x}+\bar{y})+\bar{w}(\bar{x}+\bar{z})}\tilde{J}_A(zx,\alpha(w),\alpha(y))
-(-1)^{\bar{z}(\bar{x}+\bar{y}+\bar{w})}\tilde{J}_A(zw,\alpha(x),\alpha(y))\nonumber\\
&&-(-1)^{\bar{w}(\bar{x}+\bar{y}+\bar{z})}\tilde{J}_A(xy, \alpha(z), \alpha(w))\nonumber
\end{eqnarray}
for all $x,y,z \in \mathcal{H}(A).$}
\end{Lemma}
\textbf{Proof}
The  super skew-symmetry of $\tilde{J}_A(x, y, z)$ in the homogeneous elements $w, x, y, z$ in $A$ 
follows from the super skew-symmetry of the operation " $\cdot$".
Expanding the expression in the left-hand side of (ii) and then rearranging terms, we get by the super skew-symmetry of $"\cdot"$
\begin{eqnarray}
 &&\alpha^2(w)\tilde{J}_A(x,y,z)-(-1)^{\bar{w}(\bar{x}+\bar{y}+\bar{z})}\alpha^2(x)\tilde{J}_A(y,z,w) +
(-1)^{(\bar{y}+\bar{z})(\bar{w}+\bar{x})} \alpha^2(y)\tilde{J}_A(z,w, x)\nonumber\\
&&-(-1)^{\bar{z}(\bar{x}+\bar{y}+\bar{w})}\alpha^2(z)\tilde{J}_A(w,x,y)\nonumber\\
&=&(wx\cdot\alpha(y))\cdot\alpha^2(z)+(-1)^{\bar{z}(\bar{w}+\bar{x}+\bar{y})}(\alpha(z)\cdot wx)\alpha^2(y)\nonumber\\
&&+(-1)^{(\bar{y}+\bar{z})(\bar{x}+\bar{w})}[(yz\cdot\alpha(w))\cdot\alpha^2(x)+(-1)^{\bar{x}(\bar{y}+\bar{z}+\bar{w})}(\alpha(x)\cdot yz)\cdot\alpha^2(w)]\nonumber\\
&&+(-1)^{\bar{x}(\bar{y}+\bar{z})}[(wy\cdot\alpha(z))\cdot\alpha^2(x)+(-1)^{\bar{x}(\bar{y}+\bar{z}+\bar{w})}(\alpha(x)\cdot wy)\cdot\alpha^2(z)]\nonumber\\
&&+(-1)^{\bar{z}(\bar{x}+\bar{y})+\bar{w}(\bar{x}+\bar{z})}[(zx\cdot\alpha(w))\cdot\alpha^2(y)+(-1)^{\bar{y}(\bar{x}+\bar{z}+\bar{w})}(\alpha(y)\cdot zx)\cdot\alpha^2(w)]\nonumber\\
&&-(-1)^{\bar{z}(\bar{w}+\bar{x}+\bar{y})}[(zw\cdot\alpha(x))\cdot\alpha^2(y)+(-1)^{\bar{y}(\bar{x}+\bar{z}+\bar{w})}(\alpha(y)\cdot zw)\cdot\alpha^2(x)]\nonumber\\
&&-(-1)^{\bar{w}(\bar{x}+\bar{y}+\bar{z})}[(xy\cdot\alpha(z))\cdot\alpha^2(w)+(-1)^{\bar{w}(\bar{x}+\bar{y}+\bar{z})}(\alpha(w)\cdot xy)\cdot\alpha^2(z)]\nonumber
\end{eqnarray}
Next, adding and subtracting $(-1)^{(\bar{y}+\bar{z})(\bar{x}+\bar{w})}\alpha(x)\alpha(y)\cdot\alpha(wx)$ \\(resp.
$(-1)^{(\bar{y}+\bar{z})(\bar{x}+\bar{w})}(-1)^{(\bar{y}+\bar{z})(\bar{x}+\bar{w})}\alpha(w)\alpha(x)\cdot\alpha(yz)$,\\
$(-1)^{\bar{x}(\bar{y}+\bar{z})}(-1)^{(\bar{w}+\bar{y})(\bar{x}+\bar{z})}\alpha(z)\alpha(x)\cdot\alpha(wy),$\\ 
$(-1)^{\bar{z}(\bar{x}+\bar{y})+\bar{w}(\bar{x}+\bar{z})}(-1)^{(\bar{w}+\bar{y})(\bar{z}+\bar{x})}\alpha(w)\alpha(y)\cdot\alpha(zx),$\\
$(-1)^{\bar{z}(\bar{w}+\bar{x}+\bar{y})}(-1)^{(\bar{x}+\bar{y})(\bar{z}+\bar{w})}\alpha(x)\alpha(y)\cdot\alpha(zw)$\\ and 
$(-1)^{\bar{w}(\bar{x}+\bar{y}+\bar{z})}(-1)^{(\bar{z}+\bar{w})(\bar{x}+\bar{y})}\alpha(z)\alpha(w)\cdot\alpha(xy)$) 
in the first (resp. second, third, fourth, fifth, and sixth) line of the right-hand side expression in 
the last equality above, we come to the equality (ii) of the
lemma.\\

In a Hom-Malcev superalgebra $(M=M_0\oplus M_1, \cdot,\alpha)$ we define the multilinear map $G$ by
\begin{eqnarray}
 G(w, x, y, z) &=& \tilde{J}_M (wx, \alpha(y), \alpha(z))-(-1)^{\bar{x}\bar{w}}\alpha^2(x)\tilde{J}_M (w, y, z)\nonumber\\
&&-(-1)^{\bar{w}(\bar{x}+\bar{y}+\bar{z})}\tilde{J}_M (x, y, z) \alpha^2(w)\label{s2}
\end{eqnarray}
\begin{Lemma}\label{si2} \rm{
 In a Hom-Malcev superalgebra $(M=M_0\oplus M_1, \cdot,\alpha)$ the function $G(w, x, y, z)$
defined by (\ref{s2}) is super skew-symmetric in its four variables that is
\begin{eqnarray}
 G(w, x, y, z) &=& -(-1)^{\bar{x}\bar{w}}G(x,w, y, z)\nonumber\\
 G(w, x, y, z) &=& -(-1)^{\bar{y}\bar{z}}G(w, x, z,y)\nonumber\\
 G(w, x, y, z) &=& -(-1)^{\bar{x}\bar{y}}G(w, y,x, z)\nonumber\\
 G(w, x, y, z) &=& -(-1)^{\bar{w}(\bar{x}+\bar{y})+\bar{x}\bar{y}}G(y, x, w, z)\nonumber\\
G(w, x, y, z) &=& -(-1)^{\bar{w}(\bar{x}+\bar{y}+\bar{z})+\bar{z}(\bar{x}+\bar{y})}G(z, x, y, w)\nonumber
\end{eqnarray}
for all $x,y,z,w \in \mathcal{H}(M).$}
\end{Lemma}
\textbf{Proof} Since the group $S_4=<t_{12},t_{23},t_{34}>,$ where $t_{ij}$ are transpositions between $i$ and $j,$ it suffices 
to prove $G(w, x, y, z)=-(-1)^{\bar{x}\bar{w}}G(x,w, y, z),$ $G(w, x, y, z) = -(-1)^{\bar{x}\bar{y}}G(w, y,x, z)$ and 
  $G(w, x, y, z) = -(-1)^{\bar{y}\bar{z}}G(w, x, z,y)$ for all $w,x,y,z \in \mathcal{H}(M).$\\
The first and last expression are a direct consequence of the super anticommutativity of $\cdot$ and the one of the Hom-superJacobien.\\
Consider now, the expression $G(w, x, y, z) = -(-1)^{\bar{x}\bar{y}}G(w, y,x, z)$ for all $x,y,z ,w\in \mathcal{H}(M).$ Then 
we have 
\begin{eqnarray}
 &&G(w, x, y, z)+(-1)^{\bar{x}\bar{y}}G(w, y,x, z)\nonumber\\
&=&\tilde{J}_M (wx, \alpha(y), \alpha(z))-(-1)^{\bar{x}\bar{w}}\alpha^2(x)\tilde{J}_M (w, y, z)\nonumber\\
&&-(-1)^{\bar{w}(\bar{x}+\bar{y}+\bar{z})}\tilde{J}_M (x, y, z) \alpha^2(w)+\nonumber\\
&&(-1)^{\bar{x}\bar{y}}\tilde{J}_M (wy, \alpha(x), \alpha(z))-(-1)^{\bar{x}\bar{y}+\bar{y}\bar{w}}\alpha^2(y)\tilde{J}_M (w, x, z)\nonumber\\
&&-(-1)^{\bar{x}\bar{y}+\bar{w}(\bar{x}+\bar{y}+\bar{z})}\tilde{J}_M (y, x, z) \alpha^2(w)\nonumber\\
&=&\tilde{J}_M (wx, \alpha(y), \alpha(z))-(-1)^{\bar{x}\bar{w}}\alpha^2(x)\tilde{J}_M (w, y, z)\nonumber\\
&&(-1)^{\bar{x}\bar{y}}\tilde{J}_M (wy, \alpha(x), \alpha(z))-(-1)^{\bar{x}\bar{y}+\bar{y}\bar{w}}\alpha^2(y)\tilde{J}_M (w, x, z)\nonumber\\
&=&-(-1)^{(\bar{y}+\bar{z})(\bar{w}+\bar{x})+\bar{w}\bar{x}}[ \tilde{J}_M (\alpha(y), \alpha(z),xw)+(-1)^{\bar{y}\bar{z}+\bar{x}(\bar{y}+\bar{z})}
\tilde{J}_M (\alpha(x), \alpha(z),yw)\nonumber\\
&&-(-1)^{\bar{x}\bar{w}}\tilde{J}_M(y,z,w)\alpha^2(x)-(-1)^{\bar{x}\bar{z}+\bar{y}(\bar{x}+\bar{z}+\bar{w})}\tilde{J}_M(x,z,w)\alpha^2(y)]\nonumber\\
&=&0 \mbox{  (by (\ref{c}))}\nonumber
\end{eqnarray}
The following lemma is a consequence of the definition of $G(w, x, y, z)$ and the super skew-symmetry of 
$\tilde{J}_M(t, u, v)$ and $G(w, x, y, z).$
\begin{Lemma}\label{sis3} \rm {
 Let $(M=M_0\oplus M_1, \cdot,\alpha)$ be a Hom-Malcev superalgebra. Then
\begin{eqnarray}
 &&\tilde{J}_M(wx,\alpha(y),\alpha(z))+(-1)^{\bar{w}(\bar{x}+\bar{y}+\bar{z})}\tilde{J}_M(xy,\alpha(z),\alpha(w))\nonumber\\
&&+(-1)^{(\bar{y}+\bar{z})(\bar{x}+\bar{w})}\tilde{J}_M(yz,\alpha(w),\alpha(x))
+(-1)^{\bar{z}(\bar{x}+\bar{y}+\bar{w})}\tilde{J}_M (zw,\alpha(x),\alpha(y))=0\label{s3}\\
 &&2G(w, x, y, z)-\alpha^2(w)\tilde{J}_M (x, y,z)+(-1)^{\bar{x}\bar{w}}\alpha^2(x) \tilde{J}_M(w, y, z)\nonumber\\
&&-(-1)^{(\bar{y}+\bar{z})(\bar{x}+\bar{w})}\alpha^2(y)\tilde{J}_M(z,w,x)+(-1)^{\bar{z}(\bar{x}+\bar{y}+\bar{w})}\alpha^2(z)\tilde{J}_M (w,x,y)\nonumber\\
&=&\tilde{J}_M(w x,\alpha(y),\alpha(z))+(-1)^{(\bar{x}+\bar{w})(\bar{y}+\bar{z})}\tilde{J}_M (yz,\alpha(w),\alpha(x))\label{s4}
\end{eqnarray}
for all $x,y,z,w \in \mathcal{H}(M).$}
\end{Lemma}
\textbf{Proof} From the definition of $G(w, x, y, z)$ (see (\ref{s2})) we have
\begin{eqnarray}
 &&\tilde{J}_M(wx,\alpha(y),\alpha(z))=
G(w,x,y,z)+(-1)^{\bar{x}\bar{w}}\alpha^2(x)\tilde{J}_M(w,y,z)\nonumber\\
 &&+(-1)^{\bar{w}(\bar{x}+\bar{y}+\bar{z})}\tilde{J}_M(x,y,z)\alpha^2(w),\nonumber\\
&&(-1)^{\bar{w}(\bar{x}+\bar{y}+\bar{z})} \tilde{J}_M(xy,\alpha(z),\alpha(w))=(-1)^{\bar{w}(\bar{x}+\bar{y}+\bar{z})}G(x,y,z,w)\nonumber\\
&&+(-1)^{\bar{w}(\bar{x}+\bar{y}+\bar{z})+\bar{x}\bar{y}}\alpha^2(y)\tilde{J}_M(x,z,w)
+(-1)^{(\bar{x}+\bar{w})(\bar{y}+\bar{z})}\tilde{J}_M(y,z,w)\alpha^2(x),\nonumber\\
 &&(-1)^{(\bar{x}+\bar{w})(\bar{y}+\bar{z})}\tilde{J}_M(yz,\alpha(w),\alpha(x))=(-1)^{(\bar{x}+\bar{w})(\bar{y}+\bar{z})}G(y,z,w,x)\nonumber\\
&&+(-1)^{(\bar{x}+\bar{w})(\bar{y}+\bar{z})+\bar{y}\bar{z}}\alpha^2(z)\tilde{J}_M(y,w,x)
+(-1)^{\bar{z}(\bar{x}+\bar{y}+\bar{w})}\tilde{J}_M(z,w,x)\alpha^2(y),\nonumber\\
&& (-1)^{\bar{z}(\bar{x}+\bar{y}+\bar{w})}\tilde{J}_M (zw,\alpha(x),\alpha(y))=(-1)^{\bar{z}(\bar{x}+\bar{y}+\bar{w})}G(z,w,x,y)\nonumber\\
&&+(-1)^{\bar{z}(\bar{x}+\bar{y})}\alpha^2(w)\tilde{J}_M(z,x,y)
+\tilde{J}_M(w,x,y)\alpha^2(z).\nonumber
\end{eqnarray}
Therefore, adding memberwise these four equalities and using the super skew-symmetry\\ of $"\cdot"$, $\tilde{J}_M(x, y, z)$
 and $G(w, x, y, z),$ we get (\ref{s3}).\\
Next, again from the expression of $G(w, x, y, z),$
\begin{eqnarray}
 &&\tilde{J}_M (wx, \alpha(y), \alpha(z)) +(-1)^{(\bar{x}+\bar{w})(\bar{y}+\bar{z})}\tilde{J}_M (yz, \alpha(w), \alpha(x))\nonumber\\
&=&[G(w,x,y,z)+(-1)^{\bar{x}\bar{w}}\alpha^2(x)\tilde{J}_M(w,y,z)
+(-1)^{\bar{w}(\bar{x}+\bar{y}+\bar{z})}\tilde{J}_M(x,y,z)\alpha^2(w)]\nonumber\\
&&+[(-1)^{(\bar{x}+\bar{w})(\bar{y}+\bar{z})}G(y,z,w,x)
+(-1)^{(\bar{x}+\bar{w})(\bar{y}+\bar{z})+\bar{y}\bar{z}}\alpha^2(z)\tilde{J}_M(y,w,x)\nonumber\\
&&+(-1)^{\bar{z}(\bar{x}+\bar{y}+\bar{w})}\tilde{J}_M(z,w,x)\alpha^2(y)]\nonumber\\
&=&2G(w, x, y, z)-\alpha^2(w)\tilde{J}_M (x, y, z)+(-1)^{\bar{x}\bar{w}}\alpha^2(x)\tilde{J}_M (y, z, w)\nonumber\\
&&-(-1)^{(\bar{x}+\bar{w})(\bar{y}+\bar{z})}\alpha^2(y)\tilde{J}_M (z, w, x)
+(-1)^{\bar{z}(\bar{x}+\bar{y}+\bar{w})}\alpha^2(z)\tilde{J}_M (w, x, y)\nonumber
\end{eqnarray}
so that we get (\ref{s4}).\\
\\
From Lemma and Lemma, we get the following expression of $G(w, x, y, z).$
\begin{Lemma}\label{sis4} \rm {
 Let $(M=M_0\oplus M_1, \cdot,\alpha)$ be a Hom-Malcev superalgebra. Then
\begin{eqnarray}
 G(w, x, y, z) = 2[\tilde{J}_M(wx,\alpha(y),\alpha(z))+(-1)^{(\bar{y}+\bar{z})(\bar{w}+\bar{x})}\tilde{J}_M(yz,\alpha(w),\alpha(x))]\label{s5}
\end{eqnarray}
for all $x,y,z,w \in \mathcal{H}(M).$}
\end{Lemma}
\textbf {Proof} Set
\begin{eqnarray}
&&g(w, x, y, z)=\tilde{J}_M(wx,\alpha(y),\alpha(z))+(-1)^{\bar{w}(\bar{x}+\bar{y}+\bar{z})}\tilde{J}_M(xy,\alpha(z),\alpha(w))\nonumber\\
&&+(-1)^{(\bar{y}+\bar{z})(\bar{x}+\bar{w})}\tilde{J}_M(yz,\alpha(w),\alpha(x))
+(-1)^{\bar{z}(\bar{x}+\bar{y}+\bar{w})}\tilde{J}_M (zw,\alpha(x),\alpha(y))\nonumber
\end{eqnarray}
Then (\ref{s3}) says that $g(w, x, y, z)=0$
 for all $w, x, y, z$ in $A.$ Now, by adding \\$g(w, x, y, z)-(-1)^{\bar{x}\bar{w}}g(x, w, y, z)$ to the
right-hand side of Lemma (ii), we get
\begin{eqnarray}
&& \alpha^2(w) \tilde{J}_M(x, y, z)-(-1)^{\bar{w}(\bar{x}+\bar{y}+\bar{z})}\alpha^2(x)\tilde{J}_M(y, z, w)\nonumber\\
&&+(-1)^{(\bar{y}+\bar{z})(\bar{w}+\bar{x})}\alpha^2(y)\tilde{J}_M(z, w, x)-
(-1)^{\bar{z}(\bar{x}+\bar{y}+\bar{w})}\alpha^2(z) \tilde{J}_M(w, x, y)\nonumber\\
&=&\tilde{J}_M(wx, \alpha(y), \alpha(z))+(-1)^{(\bar{y}+\bar{z})(\bar{w}+\bar{x})}\tilde{J}_M(yz, \alpha(w), \alpha(x))\nonumber\\
&&+(-1)^{\bar{x}(\bar{y}+\bar{z})}\tilde{J}_M(wy, \alpha(z), \alpha(x))
+(-1)^{\bar{z}(\bar{x}+\bar{y})+\bar{w}(\bar{x}+\bar{z})}\tilde{J}_M(zx,\alpha(w), \alpha(y))\nonumber\\
&&-(-1)^{\bar{z}(\bar{x}+\bar{y}+\bar{w})}\tilde{J}_M(zw, \alpha(x), \alpha(y))
-(-1)^{\bar{w}(\bar{x}+\bar{y}+\bar{z})}\tilde{J}_M(xy, \alpha(z), \alpha(w))\nonumber\\
&&+ \tilde{J}_M(wx, \alpha(y), \alpha(z)) + (-1)^{\bar{w}(\bar{x}+\bar{y}+\bar{z})}\tilde{J}_M(xy, \alpha(z), \alpha(w))\nonumber\\
&&+(-1)^{(\bar{y}+\bar{z})(\bar{w}+\bar{x})}\tilde{J}_M(yz, \alpha(w), \alpha(x)) + (-1)^{\bar{z}(\bar{x}+\bar{y}+\bar{w})}\tilde{J}_M(zw, \alpha(x), \alpha(y))\nonumber
\end{eqnarray}
\begin{eqnarray}
&&-(-1)^{\bar{x}\bar{w}}\tilde{J}_M(xw, \alpha(y), \alpha(z))-
(-1)^{\bar{x}(\bar{w}+\bar{y}+\bar{z})+\bar{x}\bar{w}}\tilde{J}_M(wy, \alpha(z), \alpha(x))\nonumber\\
&&-(-1)^{(\bar{y}+\bar{z})(\bar{w}+\bar{x})+\bar{x}\bar{w}}\tilde{J}_M(yz, \alpha(x), \alpha(w))
-(-1)^{\bar{z}(\bar{x}+\bar{y}+\bar{w})+\bar{x}\bar{w}}\tilde{J}_M(zx, \alpha(w), \alpha(y))\nonumber\\
&=&3\tilde{J}_M(wx, \alpha(y), \alpha(z)) + 3(-1)^{(\bar{y}+\bar{z})(\bar{w}+\bar{x})}\tilde{J}_M(yz, \alpha(w), \alpha(x))\nonumber
\end{eqnarray}
i.e
\begin{eqnarray}
&& \alpha^2(w) \tilde{J}_M(x, y, z)-(-1)^{\bar{w}(\bar{x}+\bar{y}+\bar{z})}\alpha^2(x)\tilde{J}_M(y, z, w)\nonumber\\
&&+(-1)^{(\bar{y}+\bar{z})(\bar{w}+\bar{x})}\alpha^2(y)\tilde{J}_M(z, w, x)-
(-1)^{\bar{z}(\bar{x}+\bar{y}+\bar{w})}\alpha^2(z) \tilde{J}_M(w, x, y)\nonumber\\
 &=&3\tilde{J}_M(wx, \alpha(y), \alpha(z)) + 
3(-1)^{(\bar{y}+\bar{z})(\bar{w}+\bar{x})}\tilde{J}_M(yz, \alpha(w), \alpha(x))\label{s6}
\end{eqnarray}
Next, adding (\ref{s4}) and (\ref{s6}) together, we get
\begin{eqnarray}
 &&2G(w, x, y, z)-\alpha^2(w)\tilde{J}_M (x, y,z)+(-1)^{\bar{x}\bar{w}}\alpha^2(x) \tilde{J}_M(w, y, z)\nonumber\\
&&-(-1)^{(\bar{y}+\bar{z})(\bar{x}+\bar{w})}\alpha^2(y)\tilde{J}_M(z,w,x)+
(-1)^{\bar{z}(\bar{x}+\bar{y}+\bar{w})}\alpha^2(z)\tilde{J}_M (w,x,y)\nonumber\\
&&+\alpha^2(w) \tilde{J}_M(x, y, z)-(-1)^{\bar{w}(\bar{x}+\bar{y}+\bar{z})}\alpha^2(x)\tilde{J}_M(y, z, w)\nonumber\\
&&+(-1)^{(\bar{y}+\bar{z})(\bar{w}+\bar{x})}\alpha^2(y)\tilde{J}_M(z, w, x)-
(-1)^{\bar{z}(\bar{x}+\bar{y}+\bar{w})}\alpha^2(z) \tilde{J}_M(w, x, y)\nonumber\\
&=&\tilde{J}_M(w x,\alpha(y),\alpha(z))+(-1)^{(\bar{x}+\bar{w})(\bar{y}+\bar{z})}\tilde{J}_M (yz,\alpha(w),\alpha(x))\nonumber\\
&&+3\tilde{J}_M(wx, \alpha(y), \alpha(z)) + 3(-1)^{(\bar{y}+\bar{z})(\bar{w}+\bar{x})}\tilde{J}_M(yz, \alpha(w), \alpha(x))\nonumber
\end{eqnarray}
i.e
\begin{eqnarray}
 2G(w, x,y,z)&=&4\tilde{J}_M(wx,\alpha(y),\alpha(z)) + 
4(-1)^{(\bar{y}+\bar{z})(\bar{w}+\bar{x})}\tilde{J}_M(yz,\alpha(w),\alpha(x))\nonumber\\
\mbox{ and (\ref{s5}) follows.}\nonumber
\end{eqnarray}
\section{Proof}
Relaying on the lemmas of section 2, we are now in position to prove the
theorem.\\
\\
\textbf{Proof of theorem}
First we establish the identity (\ref{s1}) in a Hom-Malcev superalgebra. We may write (\ref{s2}) in an equivalent form:
\begin{eqnarray}
\tilde{J}_M (wx, \alpha(y), \alpha(z))=(-1)^{\bar{x}\bar{w}}\alpha^2(x)\tilde{J}_M (w, y, z)+
(-1)^{\bar{w}(\bar{x}+\bar{y}+\bar{z})}\tilde{J}_M (x, y, z)\alpha^2(w)+G(w, x, y, z)\label{s7}
\end{eqnarray}
Now in (\ref{s7}), replace $G(w, x, y, z)$ with its expression from (\ref{s5}) to get
\begin{eqnarray}
&&-\tilde{J}_M (wx,\alpha(y),\alpha(z))=(-1)^{\bar{x}\bar{w}}\alpha^2(x)\tilde{J}_M (w, y, z)+
(-1)^{\bar{w}(\bar{x}+\bar{y}+\bar{z})}\tilde{J}_M (x, y, z)\alpha^2(w)\nonumber\\
&&+2(-1)^{(\bar{y}+\bar{z})(\bar{w}+\bar{x})}\tilde{J}_M (yz, \alpha(w), \alpha(x)),\nonumber
\end{eqnarray}
which leads to (\ref{s1}).
Now, we proceed to prove the equivalence of (\ref{c}) with (\ref{s1}) in a super anticommutative Hom-Malcev superalgebra. 
First assume (\ref{c}). Then Lemmas \ref{sis1}, \ref{si2}, \ref{sis3} and \ref{sis4} imply that (\ref{s1}) holds in any Hom-Malcev superalgebra.
Conversely, assume (\ref{s1}) that is
\begin{eqnarray}
 \tilde{J}_M(\alpha(y),\alpha(z),wx)&=&(-1)^{\bar{w}(\bar{y}+\bar{z})}\alpha^2(w)\tilde{J}_M(y,z,x)+
\tilde{J}_M(y,z,w)\alpha^2(x)\nonumber\\
   &&-2(-1)^{(\bar{y}+\bar{z})(\bar{x}+\bar{w})}\tilde{J}_M(\alpha(w),\alpha(x),yz)\label{as1}
\end{eqnarray}
From (\ref{as1}), we get
\begin{eqnarray}
 &&\tilde{J}_M(\alpha(y),\alpha(z),wx)+(-1)^{\bar{y}\bar{z}+\bar{w}(\bar{y}+\bar{z})}\tilde{J}_M(\alpha(w),\alpha(z),yx)\nonumber\\
&=&[(-1)^{\bar{w}(\bar{y}+\bar{z})}\alpha^2(w)\tilde{J}_M(y,z,x)+ \tilde{J}_M(y,z,w)\alpha^2(x)\nonumber\\
   &&-2(-1)^{(\bar{y}+\bar{z})(\bar{x}+\bar{w})}\tilde{J}_M(\alpha(w),\alpha(x),yz)]
+(-1)^{\bar{y}\bar{z}+\bar{w}(\bar{y}+\bar{z})}[(-1)^{\bar{y}(\bar{w}+\bar{z})}\alpha^2(y)\tilde{J}_M(w,z,x)\nonumber\\
&&+\tilde{J}_M(w,z,y)\alpha^2(x)-2(-1)^{(\bar{w}+\bar{z})(\bar{x}+\bar{y})}\tilde{J}_M(\alpha(y),\alpha(x),wz)]\nonumber\\
&=&-(-1)^{\bar{w}\bar{x}}\tilde{J}_M(y,z,x)\alpha^2(w)-(-1)^{\bar{y}(\bar{z}+\bar{w}+\bar{x})+\bar{z}\bar{w}}\tilde{J}_M(w,z,x)\alpha^2(y)\nonumber\\
&&-2(-1)^{(\bar{y}+\bar{z})(\bar{x}+\bar{w})}[\tilde{J}_M(\alpha(w),\alpha(x),yz)+
(-1)^{\bar{w}\bar{x}+\bar{y}(\bar{w}+\bar{x})}\tilde{J}_M(\alpha(y),\alpha(x),wz)]\nonumber
\end{eqnarray}
i.e
\begin{eqnarray}
&&\tilde{J}_M(\alpha(y),\alpha(z),wx)+(-1)^{\bar{y}\bar{z}+\bar{w}(\bar{y}+\bar{z})}\tilde{J}_M(\alpha(w),\alpha(z),yx)\nonumber\\
 &=&-(-1)^{\bar{w}\bar{x}}\tilde{J}_M(y,z,x)\alpha^2(w)-(-1)^{\bar{y}(\bar{z}+\bar{w}+\bar{x})+\bar{z}\bar{w}}\tilde{J}_M(w,z,x)\alpha^2(y)\nonumber\\
&&-2(-1)^{(\bar{y}+\bar{z})(\bar{x}+\bar{w})}[\tilde{J}_M(\alpha(w),\alpha(x),yz)+
(-1)^{\bar{w}\bar{x}+\bar{y}(\bar{w}+\bar{x})}\tilde{J}_M(\alpha(y),\alpha(x),wz)]\label{as2}
\end{eqnarray}
Switching simultaneous $y$ and $w$ then $x$ and $z$ in (\ref{as2}), we get
\begin{eqnarray}
&&\tilde{J}_M(\alpha(w),\alpha(x),yz)+(-1)^{\bar{w}\bar{x}+\bar{y}(\bar{w}+\bar{x})}\tilde{J}_M(\alpha(y),\alpha(x),wz)\nonumber\\
 &=&-(-1)^{\bar{y}\bar{z}}\tilde{J}_M(w,x,z)\alpha^2(y)-(-1)^{\bar{w}(\bar{x}+\bar{y}+\bar{z})+\bar{x}\bar{y}}\tilde{J}_M(y,x,z)\alpha^2(w)\nonumber\\
&&-2(-1)^{(\bar{w}+\bar{x})(\bar{z}+\bar{y})}[\tilde{J}_M(\alpha(y),\alpha(z),wx)+
(-1)^{\bar{y}\bar{z}+\bar{w}(\bar{y}+\bar{z})}\tilde{J}_M(\alpha(w),\alpha(z),yx)]\label{as3}
\end{eqnarray}
Now, replacing (\ref{as3}) in (\ref{as2}) and using the super skew-symmetry of $\tilde{J}_M,$ we get
\begin{eqnarray}
 &&\tilde{J}_M(\alpha(y),\alpha(z),wx)+(-1)^{\bar{y}\bar{z}+\bar{w}(\bar{y}+\bar{z})}\tilde{J}_M(\alpha(w),\alpha(z),yx)\nonumber\\
&=&-3(-1)^{\bar{w}\bar{x}}\tilde{J}_M(y,z,x)\alpha^2(w)-3(-1)^{\bar{y}(\bar{z}+\bar{w}+\bar{x})+\bar{z}\bar{w}}\tilde{J}_M(w,z,x)\alpha^2(y)\nonumber\\
&&+4\tilde{J}_M(\alpha(y),\alpha(z),wx)+4(-1)^{\bar{y}\bar{z}+\bar{w}(\bar{y}+\bar{z})}\tilde{J}_M(\alpha(w),\alpha(z),yx)\nonumber
\end{eqnarray}
i.e
\begin{eqnarray}
 &&-3\tilde{J}_M(\alpha(y),\alpha(z),wx)-3(-1)^{\bar{y}\bar{z}+\bar{w}(\bar{y}+\bar{z})}\tilde{J}_M(\alpha(w),\alpha(z),yx)\nonumber\\
&=&-3(-1)^{\bar{w}\bar{x}}\tilde{J}_M(y,z,x)\alpha^2(w)-3(-1)^{\bar{y}(\bar{z}+\bar{w}+\bar{x})+\bar{z}\bar{w}}\tilde{J}_M(w,z,x)\alpha^2(y)\nonumber
\end{eqnarray}
and therefore (\ref{c}) follows.

\end{document}